\documentclass[preprint,12pt]{elsarticle}
\usepackage{amssymb}
\usepackage{amsmath}
\usepackage{amsthm, algorithm}
\usepackage[noend]{algpseudocode}
\usepackage{lineno}
\usepackage{float}
\usepackage{graphicx,color}
\usepackage{wrapfig}
\usepackage{xcolor}
\usepackage{setspace}

\newcommand{\beq}{\begin{equation}}
\newcommand{\eeq}{\end{equation}}

\newcommand{\bit}{\begin{itemize}}
\newcommand{\eit}{\end{itemize}}
\newcommand{\ben}{\begin{enumerate}}
\newcommand{\een}{\end{enumerate}}

\newcommand{\by}{{\bf y}}
\newcommand{\cE}{{\mathcal{E}}} 
 
\newcommand{\bx}{ {\bf x}} 

\newcommand{\ds}{{\displaystyle}}

\def \ds          {\displaystyle}

\def \by          {{\bf y}}

\def \cg          {{\cal G}}

\graphicspath{./figures}

\title{Centroidal Voronoi Tessellation Based Methods for Optimal Rain Gauge Location Prediction}

\author[label1]{Zichao (Wendy) Di\corref{cor1}}
\ead{wendydi@mcs.anl.gov}
\author[label3]{Viviana Maggioni}
\author[label3]{Yiwen Mei}
\author[label4]{Marilyn Vazquez} 
\author[label5]{Paul Houser}
\author[label6]{Maria Emelianenko}
\address[label1]{Mathematics and Computer Science Division, Argonne National Laboratory, Lemont, IL, USA}

\address[label3]{Department of Civil, Environmental and Infrastructure Engineering, George Mason University, Fairfax, VA, USA}
\address[label4]{Mathematical Biosciences Institute, Ohio State University, Columbus, OH, USA }
\address[label5]{Department of Geography and Geoinformation Science, George Mason University, Fairfax, VA, USA}
\address[label6]{Department of Mathematical Sciences, George Mason University, Fairfax, VA, USA}

\cortext[cor1]{Corresponding author}

\begin{document}
\graphicspath{{./figures/}}
\begin{frontmatter}

\begin{abstract}
With more satellite and model precipitation data becoming available, new analytical methods are needed that can take advantage of emerging data patterns to make well informed predictions in many hydrological applications. We propose a new strategy where we extract precipitation variability patterns and use correlation map to build the resulting density map that serves as an input to centroidal Voronoi tessellation construction that optimizes placement of precipitation gauges. We provide results of numerical experiments based on the data from the Alto-Adige region in Northern Italy and Oklahoma and compare them against actual gauge locations. This method provides an automated way for choosing new gauge locations and can be generalized to include physical constraints and to tackle other types of resource allocation problems. 
\end{abstract}


\end{frontmatter}

\section{Introduction}

Precipitation is a critical variable in the water cycle, as it provides moisture for processes such as runoff, evapotranspiration, and groundwater recharge. Knowledge of the precipitation characteristics and patterns is therefore crucial for understanding land-climate interactions, for extreme event monitoring, and for water resource management. However, accurate precipitation information at fine space and time scales is difficult to obtain, as precipitation estimates from rain gauges, ground-based radars, satellite sensors, and numerical models are all affected by significant uncertainties, which can even be amplified when exposed to non-linear land surface model physics \cite{gottschalck,hazra}.

Rain gauge (or pluviometer) networks are the only direct method to measure precipitation and provide observations with high temporal resolution. As such, they are widely accepted as the benchmark for validating remotely-sensed precipitation products (e.g.,\cite{maggioni16},\cite{ana09}). However, obtaining a spatially representative precipitation field from rain gauges may require collecting a large number of observations at several locations to include different terrain, micro-climate, and vegetation variability. This translates into placing numerous gauges, which is costly in terms of maintenance and often not feasible because of location inaccessibility. Since rain gauges cannot offer spatially continuous information of precipitation \cite{kidd}, it is desirable to place them in a strategic way, accounting for changes in spatio-temporal patterns and having a methodology that is able to adapt to climate instabilities as suggested in \cite{chapter}. This work proposes an automated method to identify the optimal locations of rain gauges in order to capture the spatial variability of precipitation systems in the region and therefore provide a spatially representative rain field.

Voronoi tessellations appear in many contexts and serve a variety of purposes. They are referred to as quantizations in the electrical engineering community \cite{agrell98}, as polygons of influence or Thiessen regions in geostatistics \cite{thiessen}, as geodesic, icosahedral or hexagonal grids in climate system modeling \cite{ringler00}. The corresponding energy functional can be thought of as variation in statistical terminology, or as a cost functional in economic terms. The basic idea of this construction is to represent a large set of data by means of few representative points (generators). 

While this general concept has been widely used in a variety of contexts, the special case of a centroidal Voronoi tessellations (CVTs), which is used to denote Voronoi diagrams with the choice of generators that minimizes the energy functional, is less known due to the difficulties associated with its construction. This concept has been gaining popularity in the recent decades in many applications areas, ranging from biology and physics to finance, economics, social science, as well as hydrology (see \cite{DGJ00},\cite{Chen04} and references therein). It is particularly well studied in the context of mesh generation, clustering, quantization, imaging, reduced order modeling and partial differential equations (PDE) applications, where a number of theoretical results has been obtained attesting for its superior qualities comparing to other competing methodologies. 

The goal of this work is two-fold. First, we want to introduce the concept of CVTs for precipitation pattern analysis and point out several recently developed numerical algorithms that help constructing CVTs in continuous and discrete settings. Second, we are applying the idea of CVTs in the context of optimal placement of rain gauges, similarly to how it was previously applied in finding optimal placement of schools, post offices, and other resources \cite{DFG}. The practical aim of this study is to develop an automated strategy that would work for an arbitrary data in any geographical location. We draw attention to several modeling assumptions that are part of the algorithm presented herein and the implications of these assumptions on the algorithm performance based on several selected datasets.

We utilize the truncated-Newton (TN) method \cite{nash2000survey}, which is a large-scale nonlinear optimization algorithm, to construct all CVT solutions. The method, described briefly in Appendix (Section \ref{app:tn}), allows to considerably speed up the calculation comparing to techniques such as Lloyd method widely used in the engineering community \cite{Lloyd}, and has advantages over previously introduced modified Lloyd formulations \cite{di2012truncated}. 

As in any CVT problem aimed at finding optimal placement of resources, one needs to have a local density estimator. When trying to optimize rain gauge placement, this density should be naturally related to the measure of local precipitation variability. In this work, we have used a variability estimator based on local covariance matrix computed at the decorrelation distance. Our aim is to develop an automated strategy that could be employed in any geographical location in the future on a need basis. This method could be applied to any precipitation dataset, whether in situ, remotely, or based on model simulations.

The article is organized as follows. Section \ref{sec: methods} provides information on the CVT methodology and formulates the problem of optimal gauge placement, comparing to previously used approaches. Section~\ref{sec: datasets} gives information on the data we used in this work. The proposed algorithm is presented in Section~\ref{sec: alg}. 
We also compare predictions given by this model with existing gauge locations in Section~\ref{sec: results}. 

\section{Methods}
\label{sec: methods}
\subsection{Voronoi and centroidal Voronoi tessellations}
\label{sec:cvt}
The idea of tessellating the region, i.e. decomposing the region into sub-regions, based on the locations of rain gauges has appeared in the early works of Thiessen \cite{thiessen}.

The construction is simple. Consider a certain geographic region $ W\subset \mathbb{R}^2$. Voronoi (Thiessen) regions $\{V_i\}^k_{i=1}$ are generated by a set of points $\{\bx_i\}^k_{i=1}\in \mathbb{R}^2$  and are defined as follows:
\[
V_i(\bx_i)=\{\bx\in W: ||\bx-\bx_i||\le ||\bx-\bx_j||, j=1,\ldots, k; j\ne i\},
\]
where $\|\cdot\|$ is any distance metric. In this work, we choose it to be the standard Euclidean norm. These regions cover the entire domain $W$ and can be formed by drawing perpendicular bisectors to the segments joining consecutive generating points.  

Given a certain desired ``density'' function $\rho(\bx)$, one can compute the tessellation error, i.e. 
\beq
\ds \cE \left(\{\bx_i\}^k_{i=1}\right)=\min_{\{\bx_i\}_{i=1}^k}\sum^k_{i=1}\int_{V_i}{\rho(\bx)||\bx-\bx_i||^2 \, d\bx}.
\label{CVT_energy}
\eeq
It can be shown \cite{burkardt} that this happens precisely when $\bx_i=\bx^*_i$, where $\bx^*_i$ is the mass centroid of the corresponding Voronoi region $V_i$.
A tessellation satisfying this property is called a ``centroidal Voronoi tessellation'', or CVT for short. 
Notice that the above formulation may be extended to other more general cases, i.e. by considering other distance metrics $f(\|\cdot\|)$, other geometric constraints and periodic extensions \cite{DFG, ZED}. 

The density function may be used to represent a variety of physical characteristics, such as local characteristic length-scale \cite{Ju08}, signal intensity \cite{E10}, desired grid resolution \cite{burkardt}. In this work we propose to use it for representing spatial rainfall variability, as described in Section~\ref{sec: alg}. 

The classical method for constructing CVTs is the algorithm developed by Lloyd in the 1980s \cite{Lloyd} which represents a fixed-point type iterative mechanism. More efficient methods for calculating CVTs have been developed in the past decades \cite{DE06,DE08,di2012truncated}.  For an overview of CVT related numerical techniques we refer interested reader to \cite{DFG, Chen11}. For the purposes of this work, we will use one of the recently developed and well tested CVT solvers based on truncated Newton (TN) methodology \cite{di2012truncated} (see the preudocode given in Appendix B). 

\subsection{Problem formulation}
\label{formulation}

There are several possible ways to formulate the problem of optimal placement of rain gauges in a certain region. 

One method has been proposed in \cite{Okabe} and later used in \cite{DFG}. If $k$ rain gauges locations in region $W\in \mathbb{R}^2$ are given by $\bx_1,\ldots, \bx_k$ and $V_i\in W$ are the Voronoi (Thiessen) region associated with the $i$-th gauge, one can minimize the expected squared approximation error for the amount of rainfall $z(\bx)$ treated as a random variable $z(\bx)=m(\bx)+\epsilon(\bx)$ with mean $m(\bx)$ and deviation $\epsilon(\bx)$. If the change in the average trend $m(\bx)$ is small compared with the variance $Var(\bx)=E(\epsilon(\bx)^2)=\varepsilon$, which is assumed to be constant, the expected squared approximation error can be approximated as 
\beq
\label{opt}
\ds \cE\left(\{\bx_i\}^k_{i=1}\right) \approx \min_{\{\bx_i\}^k_{i=1}}\sum^n_{i=1}\int_{V_i} {2 \varepsilon^2 \left[1-Corr(\bx,\bx_i)\right]\, d\bx} , 
\eeq
where
\beq
\ds Corr(\bx,\bx_i)=E[\epsilon(\bx)\epsilon(\bx_i)]/[Var(\bx)Var(\bx_i)]
\label{corr}
\eeq
is the (Pearson) linear correlation coefficient of the time series at locations $\bx$ and $\bx_i$. Details on this derivation are given in Appendix A. Evaluation of \eqref{opt} is computationally intensive, requiring calculation of pairwise correlations for all points inside the domain. 

Notice that under the additional assumption that $Corr(\bx,\bx_i)$ only depends on the differences $||\bx-\bx_i||$, $F(\{\bx_i,V_i\}^k_{i=1})$ can be thought of as a generalization of the CVT energy \eqref{CVT_energy} given by:
\beq
\cE(\{\bx_i\}^k_{i=1})=\sum^k_{i=1}\int_{V_i}{\rho(\bx)f(||\bx-\bx_i||) \, d\bx} 
\label{CVT_energy1}
\eeq
with distance metric $f(||\bx-\bx_i||)=2\varepsilon^2 (1-Corr(\bx,\bx_i))$ and $\rho(\bx)=1$. The combination of conditions of constant mean and variance together with $Corr(\bx,\bx_i)=Corr(\bx-\bx_i)$ is normally referred to as the {\it weak stationarity} assumption that might or might not hold in practice, as discussed for instance in \cite{chapter}. Effectively, in this formulation the distance plays the most important role, placing a small weight on highly correlated points and magnifying weakly correlated regions. The maximum value of the variance can be rather small depending on the data, which may possibly lead to slow convergence for commonly used numerical algorithms \cite{DEJ06}.

In contrast with the above approach, standard CVT formulation \eqref{CVT_energy} applied to the same problem allows to achieve a similar effect by fixing Euclidean density and instead selecting appropriate density function $\rho$. This method is grounded on the observation that for the solutions to \eqref{CVT_energy}, the sizes of 2-dimensional Voronoi regions defined as $\ds h_{V_i}=2\max_{\by\in V_i}{||\bx_i-\by||}$ satisfy \cite{DFG}:
\beq
\ds
\frac{h_{V_i}}{h_{V_j}}=\left( \frac{\rho(\bx_i)}{\rho(\bx_j)}\right)^{1/3}.
\label{ratio}
\eeq
This is the approach advocated for in this work. One may rescale the density $\rho$ to achieve any desired ratio of cell sizes based on a certain spatial distribution of interest, for instance, in \cite{Ju08}, CVT mesh was generated using a velocity field. In this work we choose our density so that it satisfies the following conditions: \sloppy $\ds \sup_{\bx \in W}{\rho(\bx)}= \lim_{Corr(\bx)\rightarrow C_{min}} \rho(\bx) = R$ and $\ds \inf_{\bx \in W}{\rho(\bx)}= \lim_{Corr(\bx)\rightarrow C_{max}} \rho(\bx) = r$.
Here $R$ and $r$ are scale parameters and $Corr(\bx)$ denotes effective (average) correlation 
at spatial location $\bx$, with $-1 \le C_{min}\le Corr(\bx) \le C_{max}\le 1, \forall \bx\in W$. The method for computing this quantity based on averaging correlations at a ``decorrelation'' distance is discussed in \ref{sec: eff}. There are many functional forms such a relation can take. One choice is to consider power law relation of the type 
\begin{equation}
\rho(\bx)=r+R\Bigl(\frac{C_{max}-Corr(\bx)}{C_{max}-C_{min}}\Bigr)^\alpha,
\label{rho}
\end{equation}
where $\alpha>0$ denotes the power exponent.  

This approach results in the following alternative formulation of the optimal rain gauge placement problem:

\beq
\ds
\min_{\{\bx_i\}_{i=1}^k}\sum^k_{i=1}\int_{V_i}\Big[ r+ R\Big( \frac{C_{max}-Corr(\bx)}{C_{max}-C_{min}}\Big)^\alpha\Big]\|\bx-\bx_i\|_2^2 d\bx.
\label{opt-fin}
\eeq

\begin{figure}[t]
\includegraphics[width=0.5\textwidth]{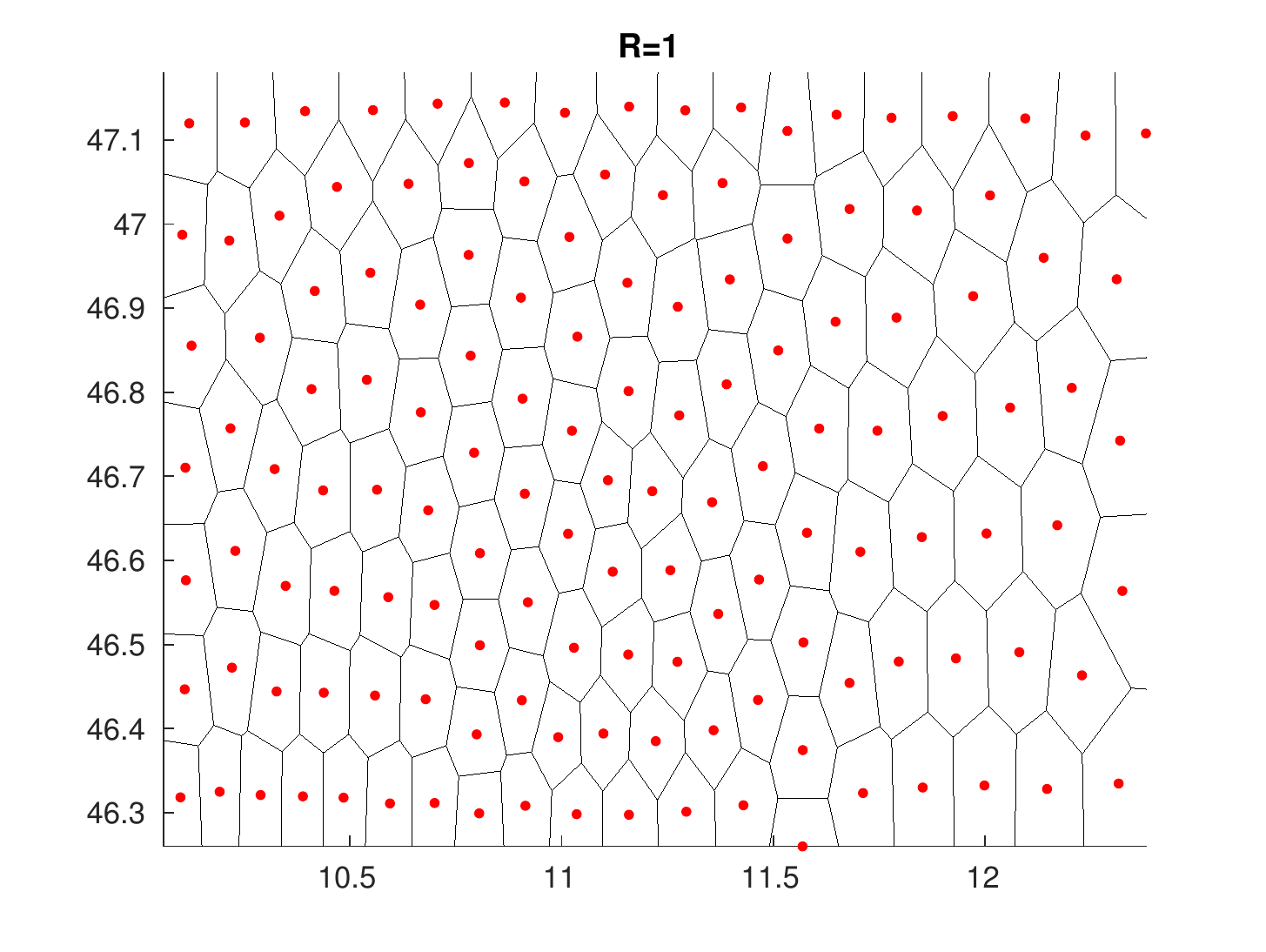}
\includegraphics[width=0.5\textwidth]{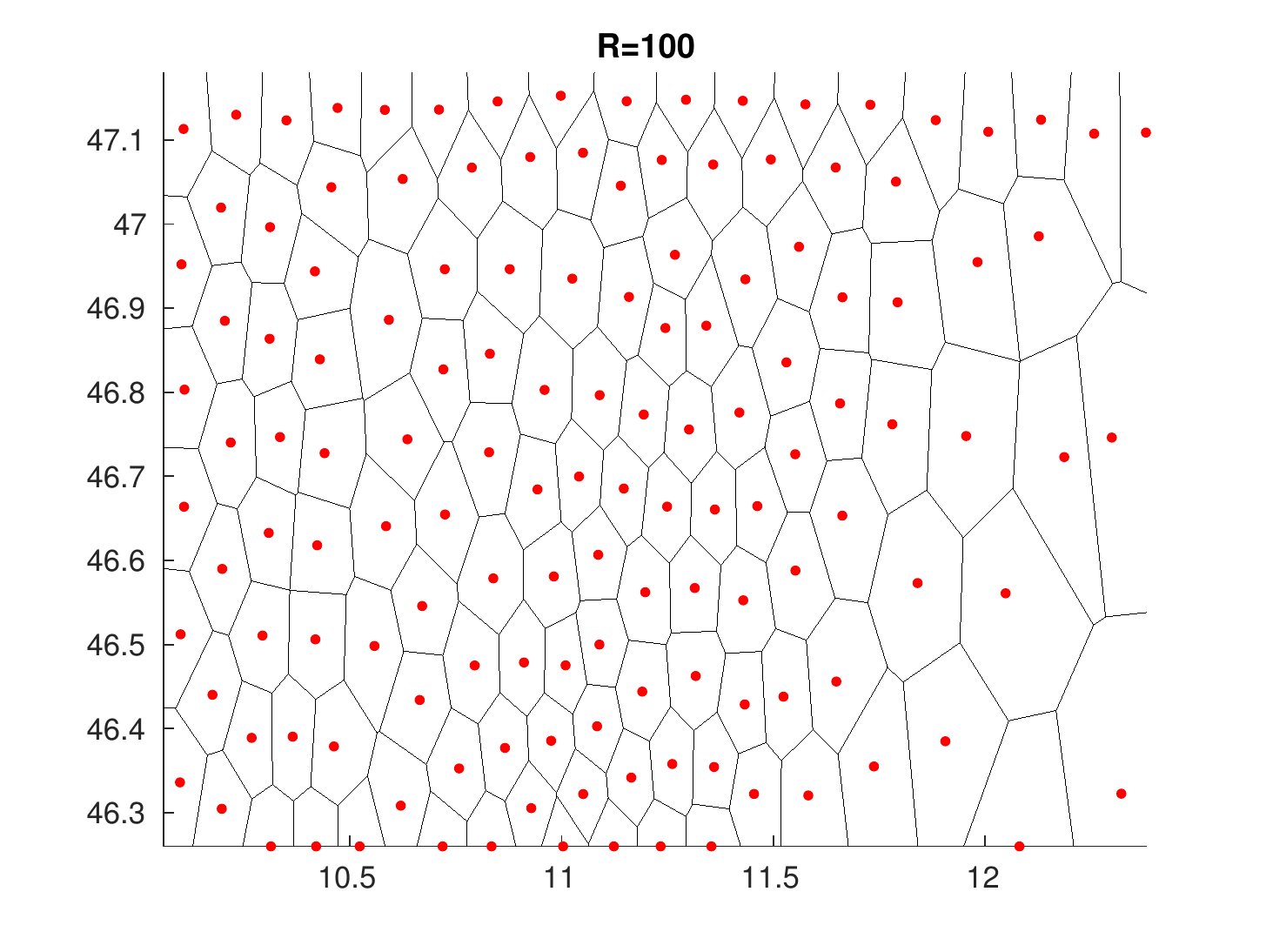}
\caption{Illustration of the effect of the ratio of $R/r$ on the sizes of Voronoi regions. $r=1$ is fixed for both figures. }
\label{fig: ratio}
\end{figure}

The choice of parameters $r,R,\alpha$ is ultimately dependent on the data, the choice of the numerical optimization algorithm and particular application. Due to possible convergence issues, one needs to constrain density away from zero, so that $r>0$. The scale parameter $R$ magnifies the range of density values. Based on \eqref{ratio}, the ratio $R/r$ can be interpreted as an approximation for the ratio of largest and smallest Voronoi regions: $\ds\frac{R}{r}\approx \frac{h_{max}}{h_{min}}$. We illustrate the role of this ratio in Figure~\ref{fig: ratio} by fixing $r=1$ and varying $R$ values. We defer to the values of $r=10^{-6}$ and $R=1$ in this work to provide the proof of concept. This choice gives a ratio of $R/r=10^6$ and was optimized for the use of TN CVT construction algorithm.

Given desired region size ratio, parameter $\alpha$ allows to enhance the contrast between peaks and valleys of the function. It is expected that higher values of $\alpha$ will essentially penalize low correlation areas compared to high correlation (low density) areas, exaggerating density differences inside the given domain. The choice of the enhancement parameter $\alpha$ is described in details in Section~\ref{sec: alg}.

\section {Study Regions and Dataset}
\label{sec: datasets}

This work focused on two very different regions: Oklahoma in the United States and Alto-Adige in Northern Italy. The reason we chose these two domains is twofold. First, they are both covered by dense rain gauge networks that can be used as reference to evaluate the proposed algorithm. Second, they are characterized by different topography and therefore different precipitation processes.

\begin{wrapfigure}{r}{0.5\textwidth}
\includegraphics[width=0.5\textwidth]{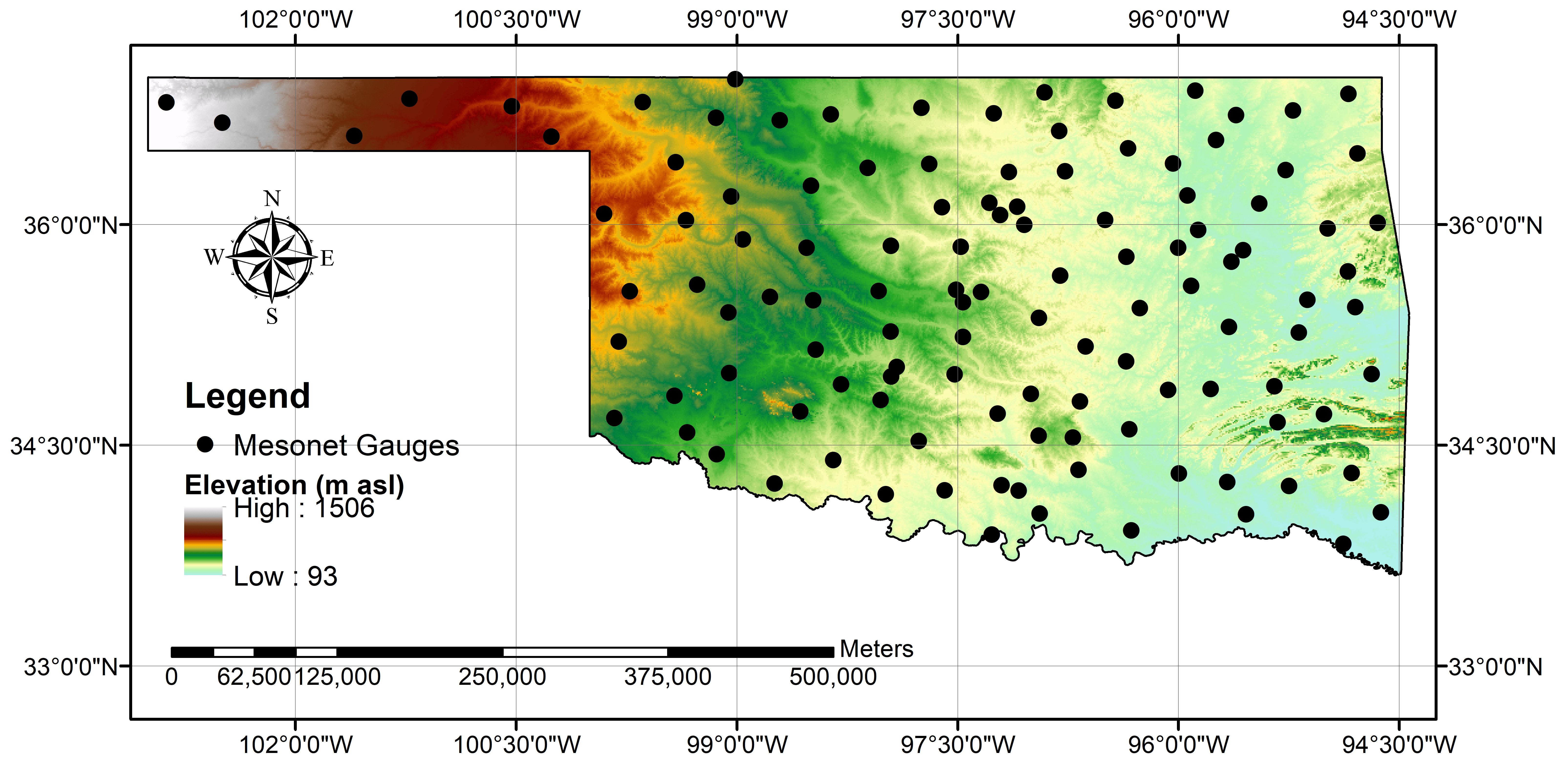}
\includegraphics[width=0.5\textwidth]{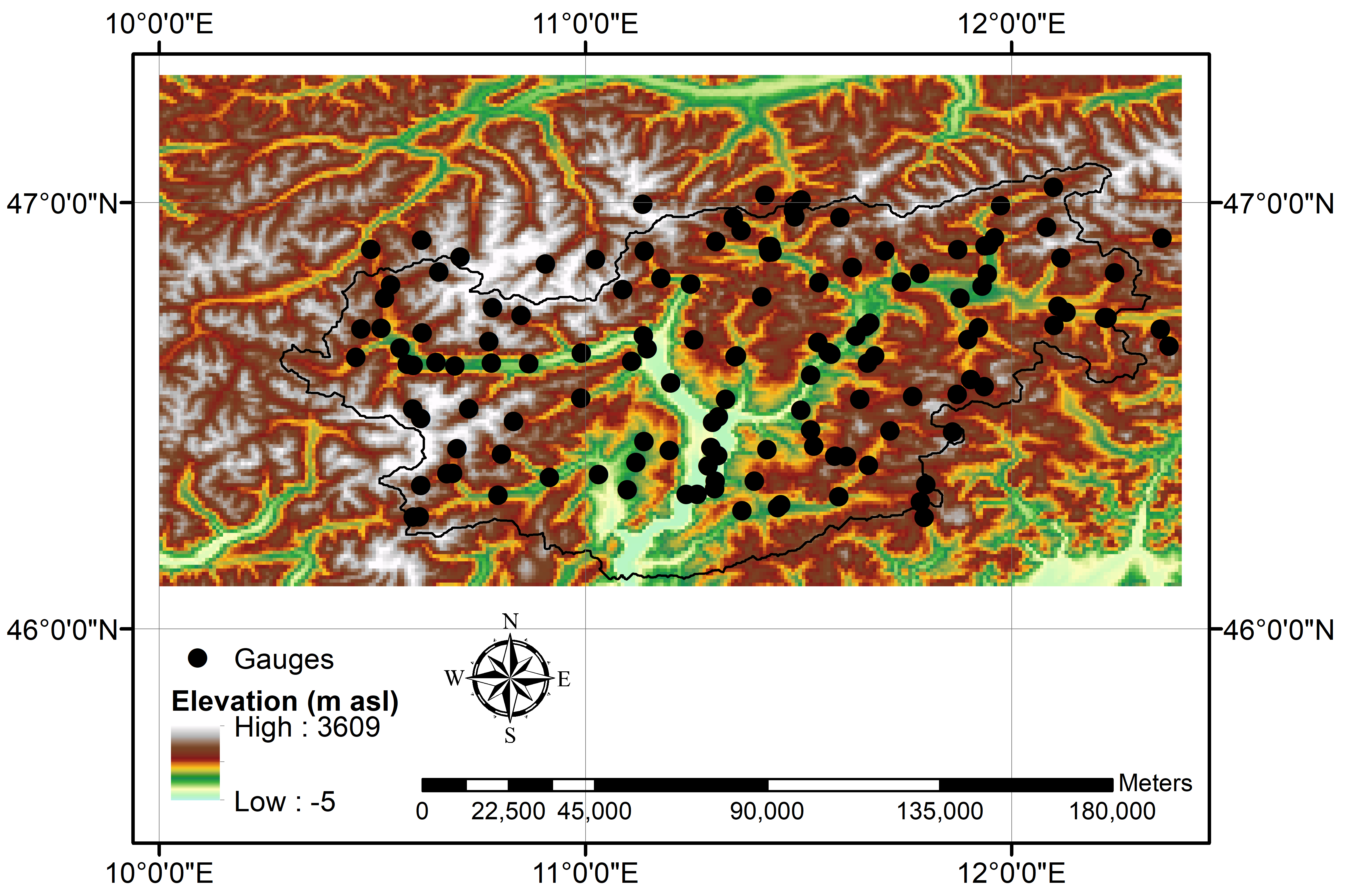}
\caption{Topographic maps with existing gauge locations: (top) Oklahoma region; (bottom) Adige region.}
\label{fig: gauges}
\end{wrapfigure}

Oklahoma is characterized by relatively uniform terrain, with gentle topography that rises from the southeastern corner to the tip of the panhandle. The continental climate of the region presents cold winters and hot summer seasons and a rainfall spatial pattern that exhibits a west-to-east (dry-to-wet) gradient. The study domain is covered by a dense network of meteorological stations, the Oklahoma Mesoscale Network (hereinafter Mesonet; \cite{brock95}), as seen in Figure~\ref{fig: gauges}~(top). Although data collected at these stations were not used in this study, their locations were compared to the output from our proposed algorithm that optimizes gauge placement to fully capture the precipitation variability in the region.

A similar analysis was performed in Alto-Adige, located in the eastern Italian Alps. Unlike Oklahoma, this area is characterized by complex topography, with elevation ranging from 65 to almost 4000 m a.s.l.

Precipitation climatology in the area exhibits strong spatial gradients, with mean annual precipitation varying from $\sim$500 mm in the Northwestern region to $\sim$1700 mm in the Southeastern part (\cite{maggioni}, \cite{niko}). A network of 192 rain gauges is available in the region. The rain gauges are distributed quite uniformly over the area (Fig.~\ref{fig: gauges}~(bottom)), providing a very dense gauge density ($\sim$1/70 rain gauge/km$^2$) for a mountainous area. For the regional studies across Oklahoma and Alto-Adige, we adopted a high resolution (1hour/0.04$^{\circ}$) satellite precipitation product, the Precipitation Estimation from Remotely Sensed Imagery using Artificial Neural Networks (PERSIANN; \cite{hsu97},\cite{sorooshian},\cite{hsu08}), produced by the University of California, Irvine. PERSIANN uses an adaptive algorithm to extract information from Infra-Red (IR) brightness temperature images, classify the extracted feature, and apply multivariate mapping of the classified texture to the rain rate. Specifically, in this study we used the PERSIANN-Cloud Classification System (PERSIANN-CCS), which only applies information from IR observations from geostationary satellites with high sampling frequency (\cite{hsu10}). In the numerical examples below, we used a crude conversion of $0.04^{\circ} = 5$[km] when presenting the results for both regions of interest. Improving this approximation or changing it for other geographic locations will bear no significant consequences in terms of the results.

\section{Proposed algorithm}
\label{sec: alg}

The main steps of the algorithm we are proposing are the following:
\ben
\item[Step 1] Build CVT density function based on the information on spatio-temporal correlations in a certain region using \eqref{rho},
\item[Step 2] Solve the CVT minimization problem \eqref{opt-fin} to obtain optimal rain gauge locations for this region. 
\een

We are now going to discuss the details of Step 1, while an example of a TN-based CVT solver is provided in Appendix B. Note that other methods may be used in place of the TN method and might be equally effective depending on the properties of the data.

\begin{figure}[h]
    \centering
    \includegraphics[width=0.51\textwidth]{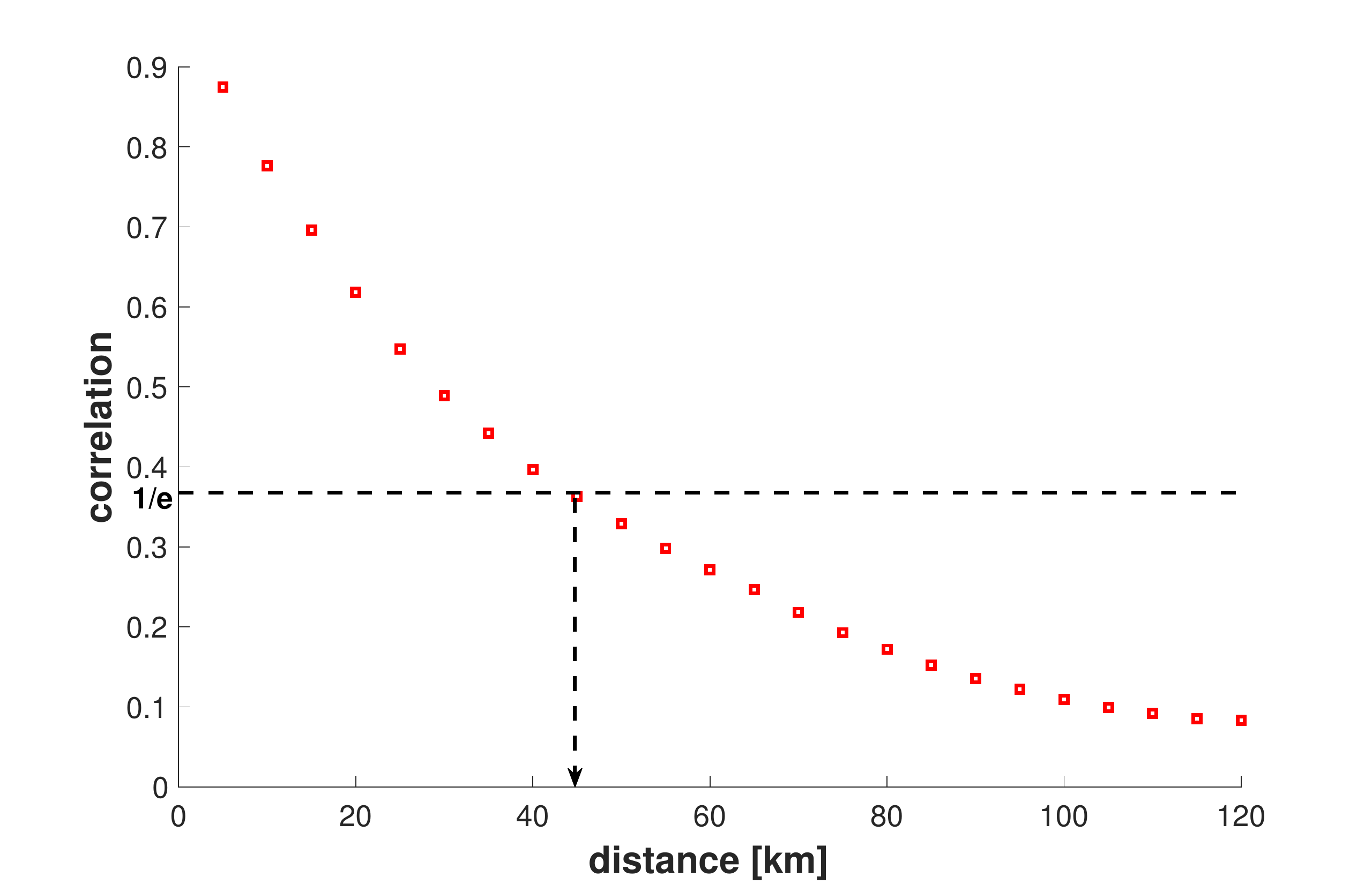}
    \includegraphics[width=0.45\textwidth]{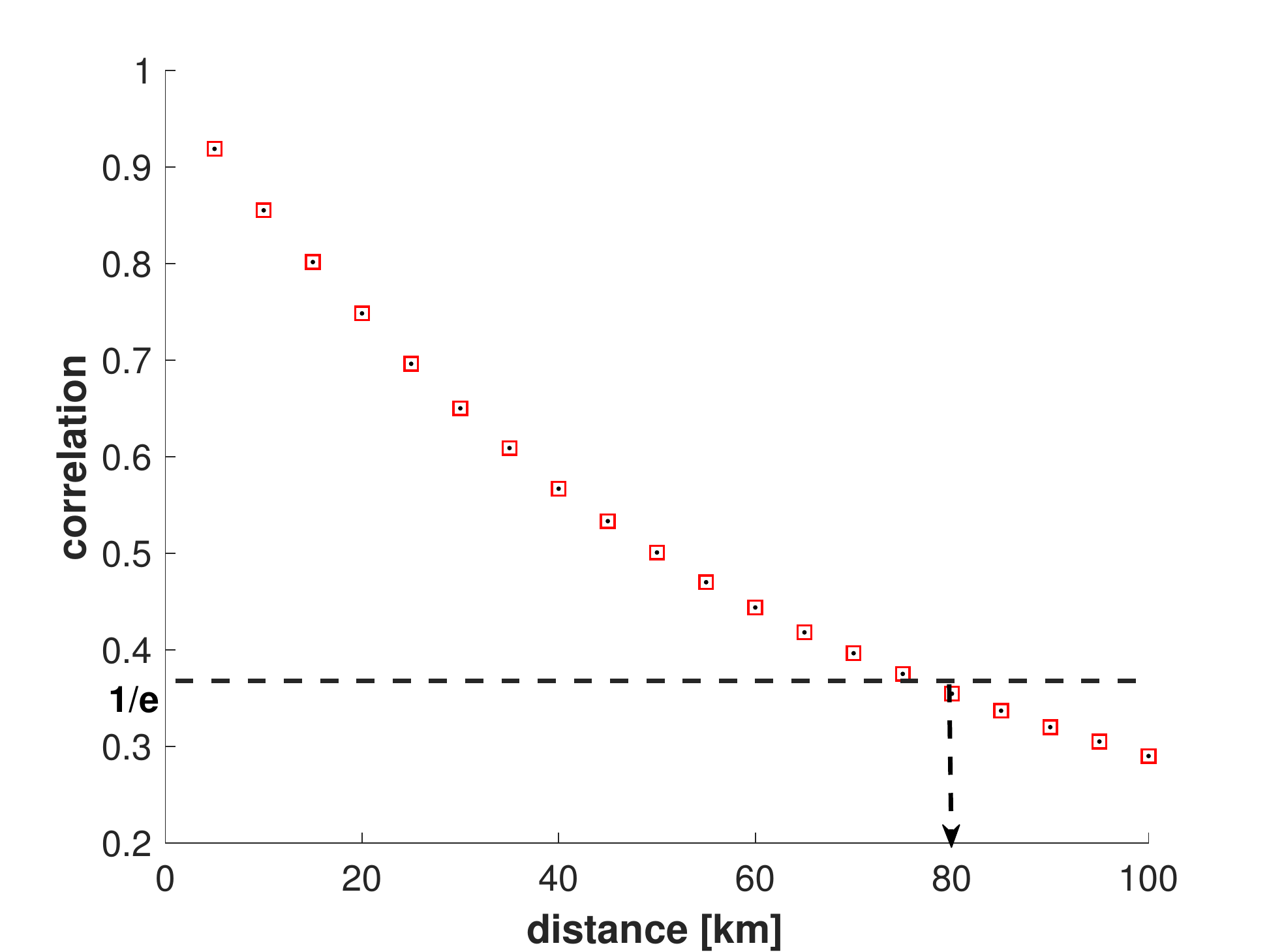}
    \caption{Calculation of the decorrelation distance using $1/e$-rule: (left) for the mountainous region in Northern Italy, with $d=45$ [km]; (right) for the Oklahoma region, with $d=80$ [km].}
    \label{fig:decor}
\end{figure}

 \subsection{Effective correlation computation}
\label{sec: eff}

 Effective local correlation $Corr(\bx)$ is a key ingredient of the CVT density estimator used in our method. The approach we take is based on the calculation of local correlation of the time series at each grid point. The number of neighbors one should be taking into account in this computation is related to the {\it decorrelation distance (radius)}, which is estimated using the the exponential model with the so-called nugget effect \cite{ciach99,ciach06}:
\begin{equation}
\rho_g(d)=c_0 \exp{\big[ -\big( \frac{d}{d_0}\big)^{s_0}\big]},
\label{rhog}
\end{equation}
where $d$ is the separation distance between two points, $c_0$ is the nugget parameter (which corresponds to the correlation value for the near-zero distances; \cite{cressie}), $d_0$ is the scale parameter (which corresponds to the spatial decorrelation distance), and $s_0$ is the correlogram shape parameter, which controls the behavior of the model near the origin for small separation distances. Note that $(1-c_0)$ is the instant decorrelation due to random errors in the rainfall observations \cite{ciach03}.

 We estimate $Corr(\bx)$ as the average of correlation coefficients between given point $\bx$ and locations on the circle $S(\bx,d)=\{\by\in W : ||\by-\bx||_2 = d \}$ of radius $d$:

\beq
Corr(\bx) = Corr_d(\bx)=\frac{1}{| S(\bx,d)|}\sum_{\by\in S(\bx,d)}{Corr(\bx,\by)},
\eeq
where $Corr(\bx,\by)$ is the Pearson correlation coefficient given in \eqref{corr}.

Monte Carlo integration over the circle of radius $d$ may be used to speed up the calculation of this quantity. Using uniform sampling $\by_1,\ldots, \by_N$ over the region $\tilde S(\bx,d)=\{\by\in W : d-1\le ||\by-\bx||_2\le d+1 \}$, we define

\beq
Corr_N(\bx) \approx \frac{1}{N}\sum^N_{i=1}{Corr(\bx,\by_i)}
\label{corr_N}
\eeq

This method gives an error of the order of $O(1/\sqrt{N})$ \cite{caflisch_1998}. The value of $N=100$ was used in the numerical experiments described in Section~\ref{sec: results}. Importance sampling may be used to improve the accuracy of the Monte Carlo approximation of the local correlation map, but it was not explored in current work. 

We define the separating distance $d_0$ at which the correlation is $1/e$ the correlation length for the (assumed) exponential variogram model:
\beq
d_0 = \min\{d: Corr_d(\bx)<\frac{1}{e} \}
\label{decorr}
\eeq

\begin{figure}[h]
\centering
    \includegraphics[width=0.9\linewidth]{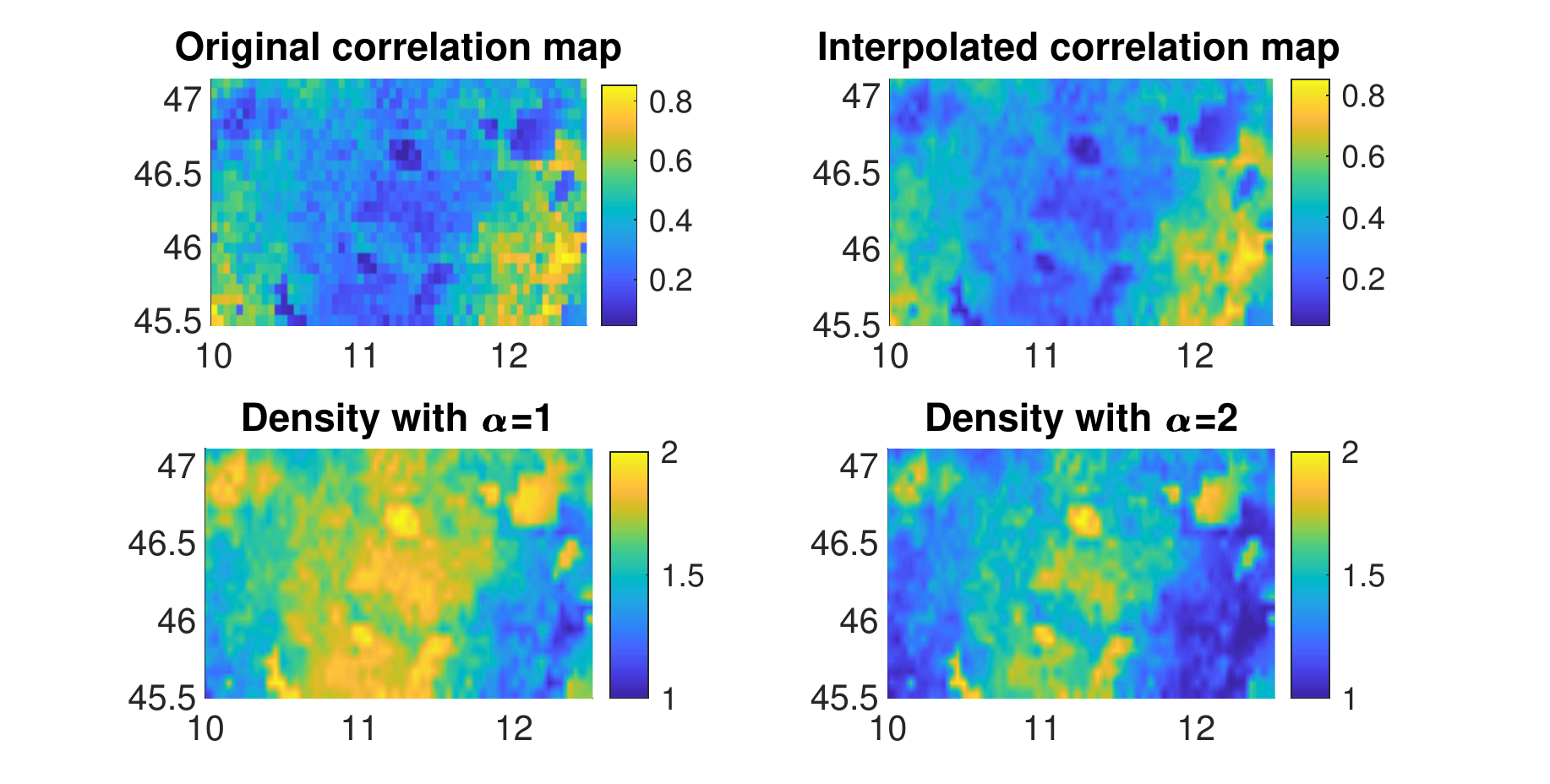}
    
    \caption{Density computation for the Italy region: (top left) Effective correlation map based on the original data; (top right) interpolated map; (bottom left) density map constructed using \eqref{rho} with $\alpha=1$; (bottom right) density map for $\alpha=2$.}
    \label{fig: Adige-dens}
\end{figure}


 Figure~\ref{fig:decor} (left) demonstrates the calculation of the decorrelation distance for the case of the data collected over the Adige mountainous region in Northern Italy. We expect to have relatively small decorrelation distance over mountains, which is indeed true based on this data, from which a value of $9$ grid points may be estimated, which corresponds to 45 [km]. Similar calculation performed over the flat region of Oklahoma shows decorrelation at $16$ grid points, corresponding to 80 [km], as seen in Figure~\ref{fig:decor} (right).

Now that the correlation component of the density is defined, we are ready to discuss a possible strategy for picking optimal value for $\alpha$ given expected number of rain gauges and a correlation threshold.

The choice of the density can be motivated by several factors, including availability of resources and desired accuracy.

Figures~\ref{fig: Adige-dens} 
shows CVT densities for Adige region, computed using formula \eqref{rho} with $Corr(\bx)$ computed at the decorrelation distance as discussed above. Two different values of parameter $\alpha$ have been tested, resulting in clear differences in densities.




\subsection{Optimal $\alpha$ and main algorithm}

One may want to optimize the choice of parameter $\alpha$ based, for instance, on the desired number of rain gauges $k_g$ to be placed in the regions with relative correlation below certain threshold $C_{tol}$. Namely, we can pick density in such a way that we get approximately $k_g$ locations with relative correlation below $C_{tol}$.  We may choose optimal value $\alpha^*$ as the smallest $\alpha$ satisfying

\begin{equation}
k=\Big|\{\bx\,|\,C^\alpha_{rel}=\Big(\frac{Corr(\bx)-C_{min}}{C_{max}-C_{min}}\Big)^\alpha < C_{tol}\} \Big| \le k_g\label{alpha_opt}   
\end{equation}
where $|\cdot|$ denotes the number of grid points in the set. 

If $\alpha=1$ immediately satisfies condition \eqref{alpha_opt}, we stop. Otherwise we keep increasing $\alpha$ until desired resolution is obtained. 

For instance, if the desired relative tolerance is $C_{tol}=0.1$, meaning that locations with correlation below 10\% are targeted, and if we are planning to place $k_g=100$ gauges in a certain area, we will need to pick $\alpha$ to satisfy
$$
k=\Big|\{\bx\,|\,C^\alpha_{rel} < 0.1\Big| \le 100
$$

\begin{figure}[t]
    \includegraphics[height=3cm, trim={0.5cm 0.5cm  0.5cm 0.5cm},clip]{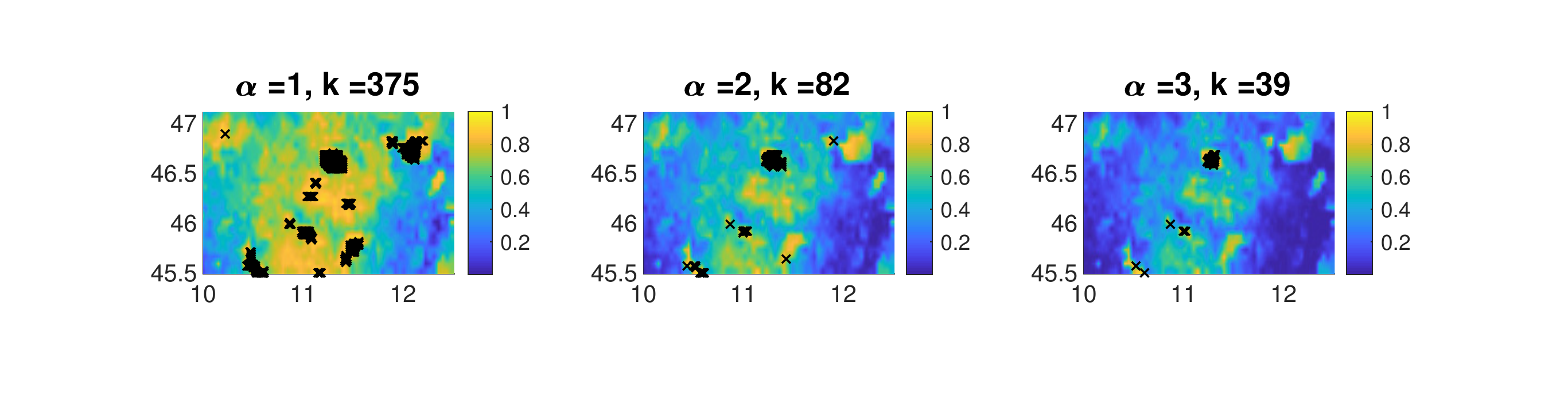}
    \caption{Comparing density computation for different $\alpha$ values: Adige region data. Points satisfying condition \eqref{alpha_opt} are marked and their number is denoted as $k$ in the title to each subfigure.}
    \label{fig: Adige3figs}
    \end{figure}

Figure~\ref{fig: Adige3figs} gives a visualization of the optimal $\alpha$ selection procedure in the case of Adige data based on the above strategy, with the choice of $k_g=60, C_{tol}=0.1$ and default values $r=1,R=1$. 

While the locations satisfying condition \eqref{alpha_opt} are good candidates for initial gauge placement, this choice is not optimal in terms of the overall approximation of the precipitation in the region of interest. As discussed above, optimal placement is attained by computing CVT that minimizes the approximation error \eqref{CVT_energy}.

Finally, in Algorithm \ref{alg} we describe the main iterative algorithm for determining optimal placement for rain gauges for a given region.

\begin{algorithm}[h]
\caption{Automatic selection of optimal gauge locations.}
\label{alg}
\begin{algorithmic}[1]
\Procedure{GaugeOptim}{}
\State Define correlation threshold $C_{tol}$ (default value $C_{tol}=0.1$), desired number of gauges $k_{g}$.
\State 
Form the $N\times n$ observation matrix $Y$, where $N$ and $n$ denote the discrete spatial and temporal resolution, respectively. 

\State Compute $Corr(\by_j)$ at each location $\by_j$, $j=1,\ldots N$ using \eqref{corr_N}. 
\State Compute decorrelation distance $d$ using \eqref{decorr}.
\State Interpolate the correlation map. Set $\alpha=0$.
\State \emph{loop}: Let $\alpha=\alpha+1$.
\State Build the density map for the interpolated grid using \eqref{rho}.
\If {\eqref{alpha_opt} is not satisfied} 
\State \textbf{goto} \emph{loop}.  
\Else
\State Return optimal $\alpha$.
\EndIf

\State  Starting with $k_g$ random points, find CVT-optimal generators $\bx_i, i=1,\ldots,k_{g}$ using TN method given in Appendix B or any alternative method with the density \eqref{rho}.
\State Calculate CVT energy \eqref{CVT_energy}. 
\EndProcedure
\end{algorithmic}
\end{algorithm}

Apart from the precipitation time series data, 
the only input parameters needed to run the code are the correlation threshold $C_{tol}$ and the number of gauges to be placed $k_{g}$. The default value for $C_{tol}$ is 0.1, while the number of gauges can be arbitrary. If desired, the user may choose to construct the CVT tessellation for any given number of generators $k_g$ starting immediately at line 13 of Algorithm~\ref{alg}.

\section{Numerical results and discussion}
\label{sec: results}

\begin{figure*}[h]
\includegraphics[trim={0 0.1cm  0 0},clip]{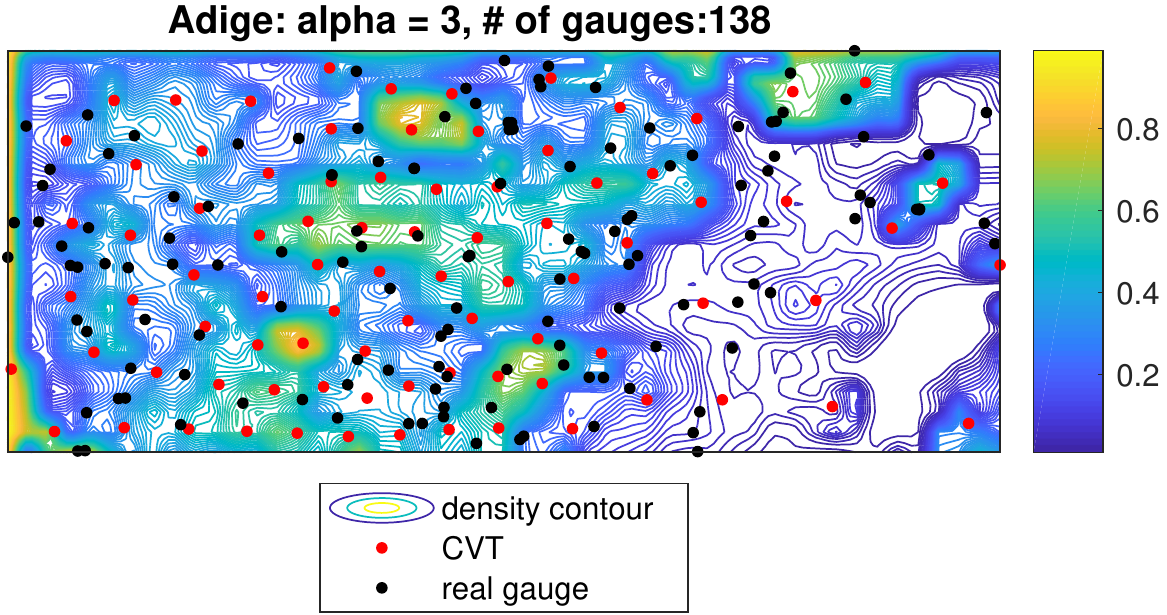}
\caption{Optimized rain gauge locations compared with real gauge data for the Adige dataset. }
\label{fig:real1}
\end{figure*}

\begin{figure*}[h]
\includegraphics[trim={0.2cm 1.6cm  0.2cm 0.4cm},clip]{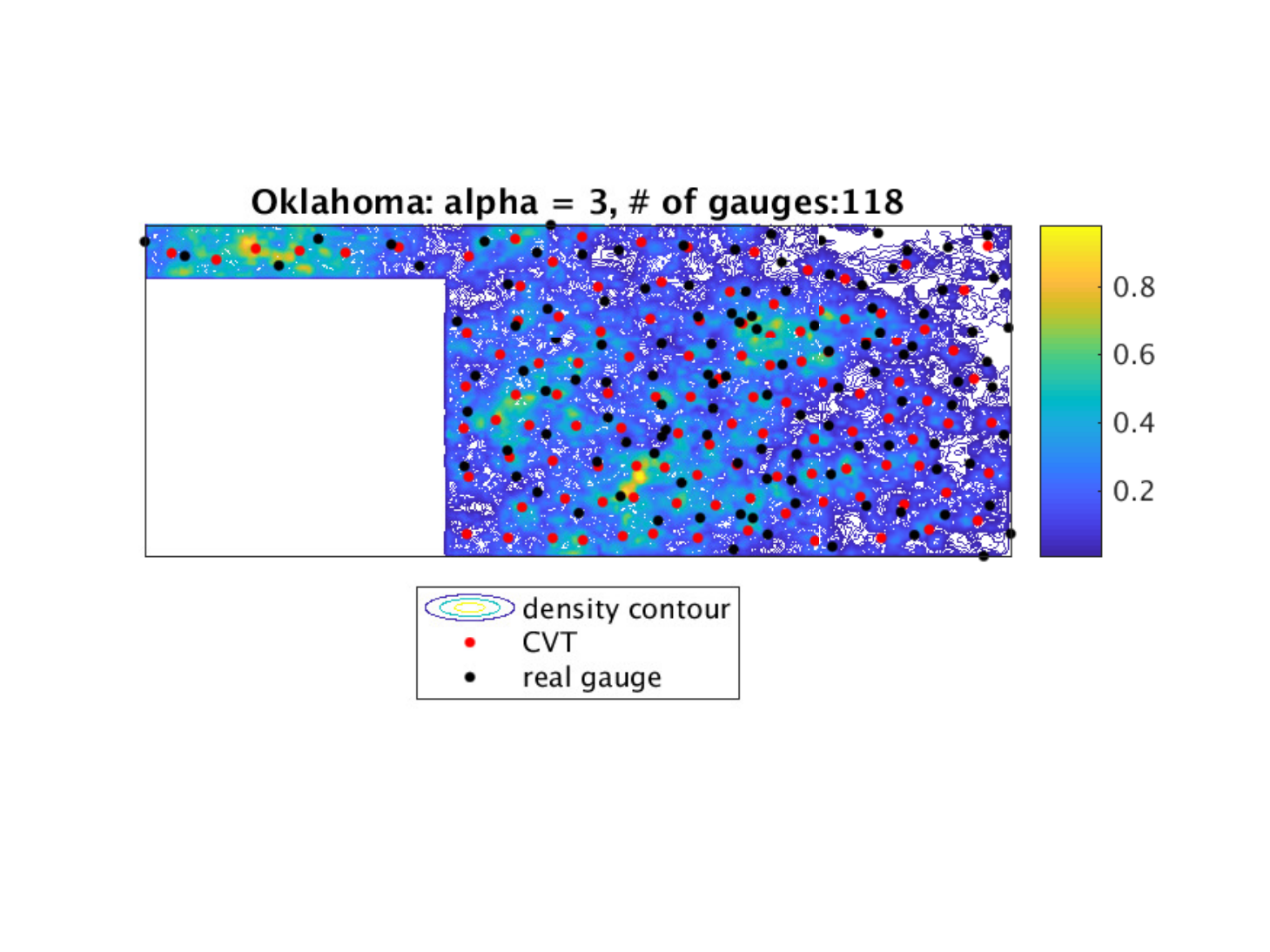}
\caption{Optimized rain gauge locations compared with real gauge data for the Oklahoma dataset. }
\label{fig:real2}
\end{figure*}

Figures~\ref{fig:real1} and \ref{fig:real2} provide the results of applying Algorithm~\ref{alg} to both regions and compare the optimal locations of rain gauges with existing rain gauge locations. The default values of $r=10^{-6},R=1$ were chosen in all calculations. 

As shown in Figure ~\ref{fig:real1}, more gauges would be required where topography is the most complex, i.e., the northwestern region of Alto Adige, characterized by altitudes close to 4000 m a.s.l. On the other hand, in the valleys in the southeastern region, less gauges would be needed to fully capture the variability of precipitation. Similarly, in Oklahoma, where topography is more uniform than in Alto Adige, the gauge placement is also more homogeneous, with limited areas characterized by high density. This confirms that the proposed algorithm is able to optimize gauge placement based on observed precipitation patterns, which are linked to the geography of the region (e.g., orographic rainfall systems). In fact, elevation was shown in past literature to explain most of the variance in rainfall in mountainous regions (e.g., \cite{sanchez, johnson}).

\begin{figure*}
\includegraphics[width=0.5\textwidth]{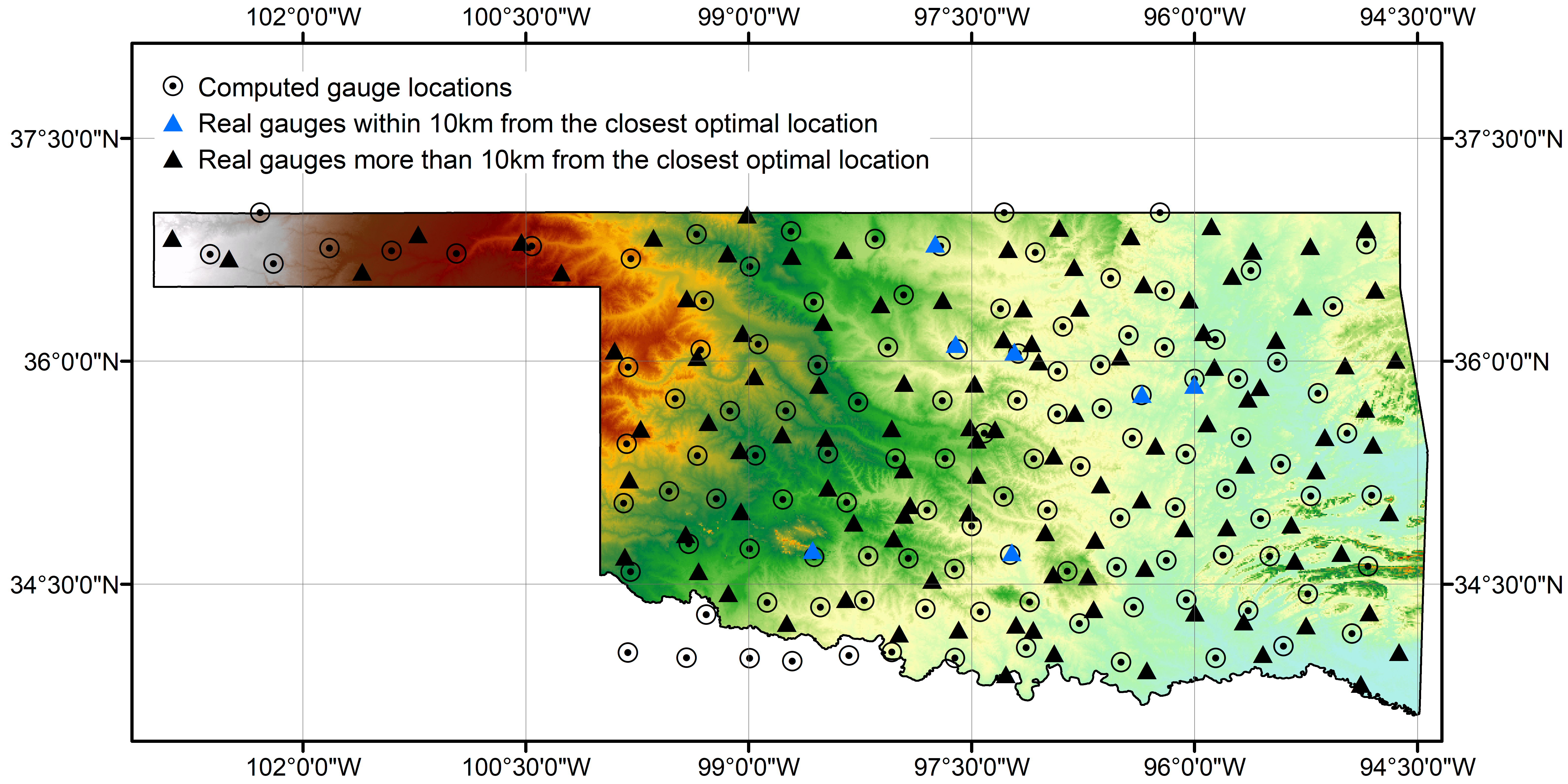}
\includegraphics[width=0.5\textwidth]{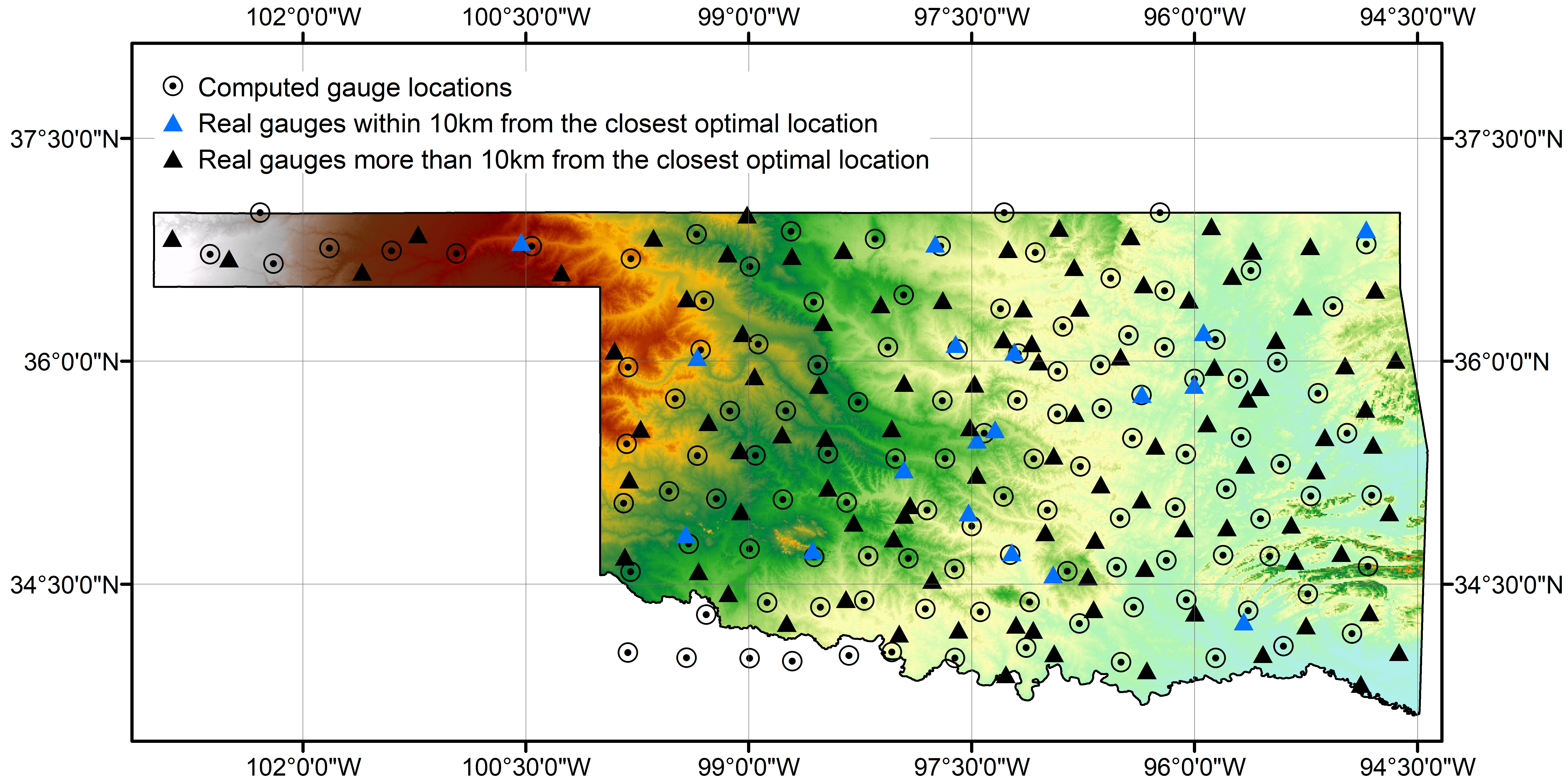}
\caption{Comparing existing locations with optimal placement for the Oklahoma region. }
\label{fig:compare_OK}
\end{figure*}

\begin{figure*}[ht]
\includegraphics[width=0.5\textwidth]{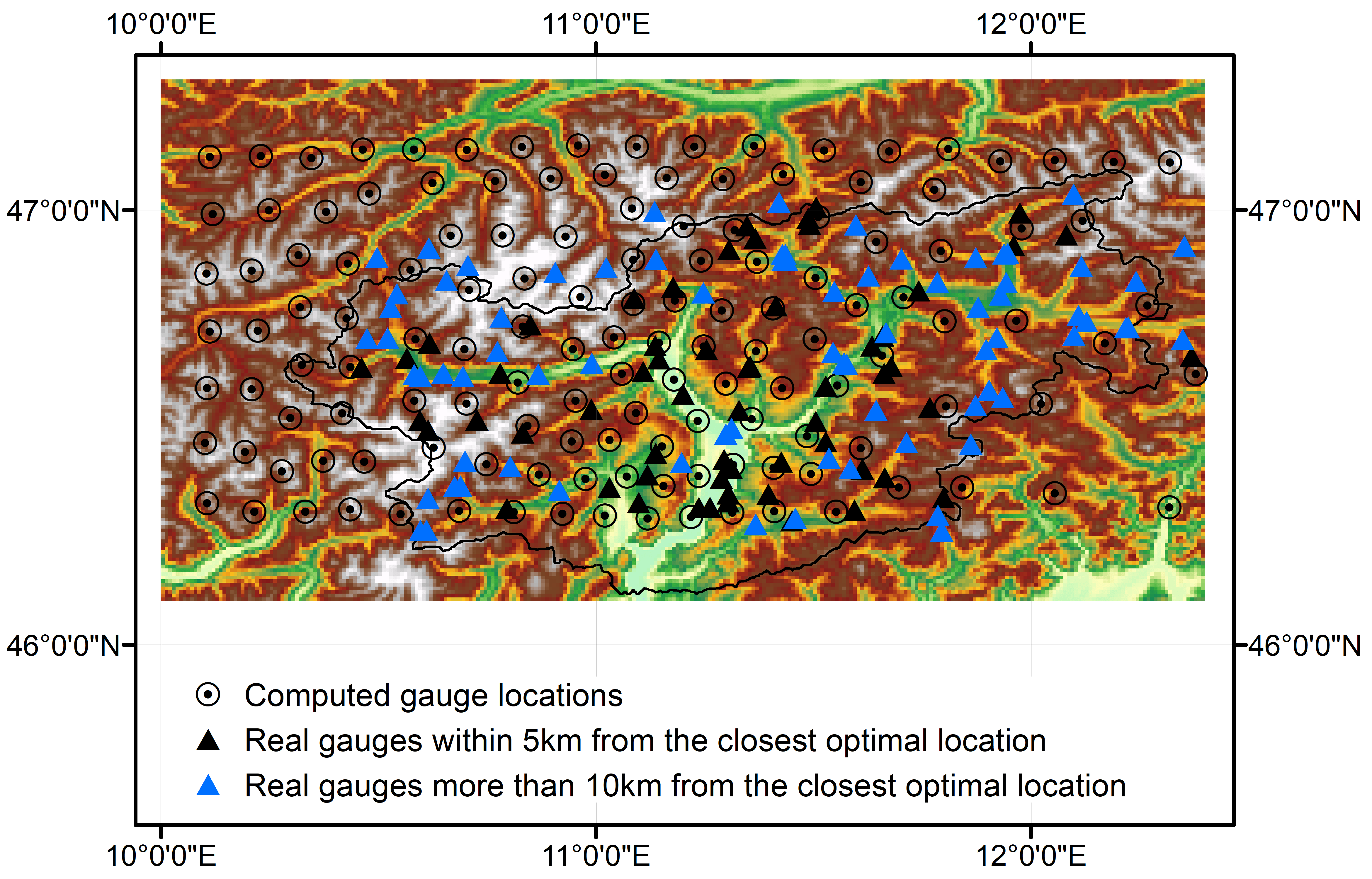}
\includegraphics[width=0.5\textwidth]{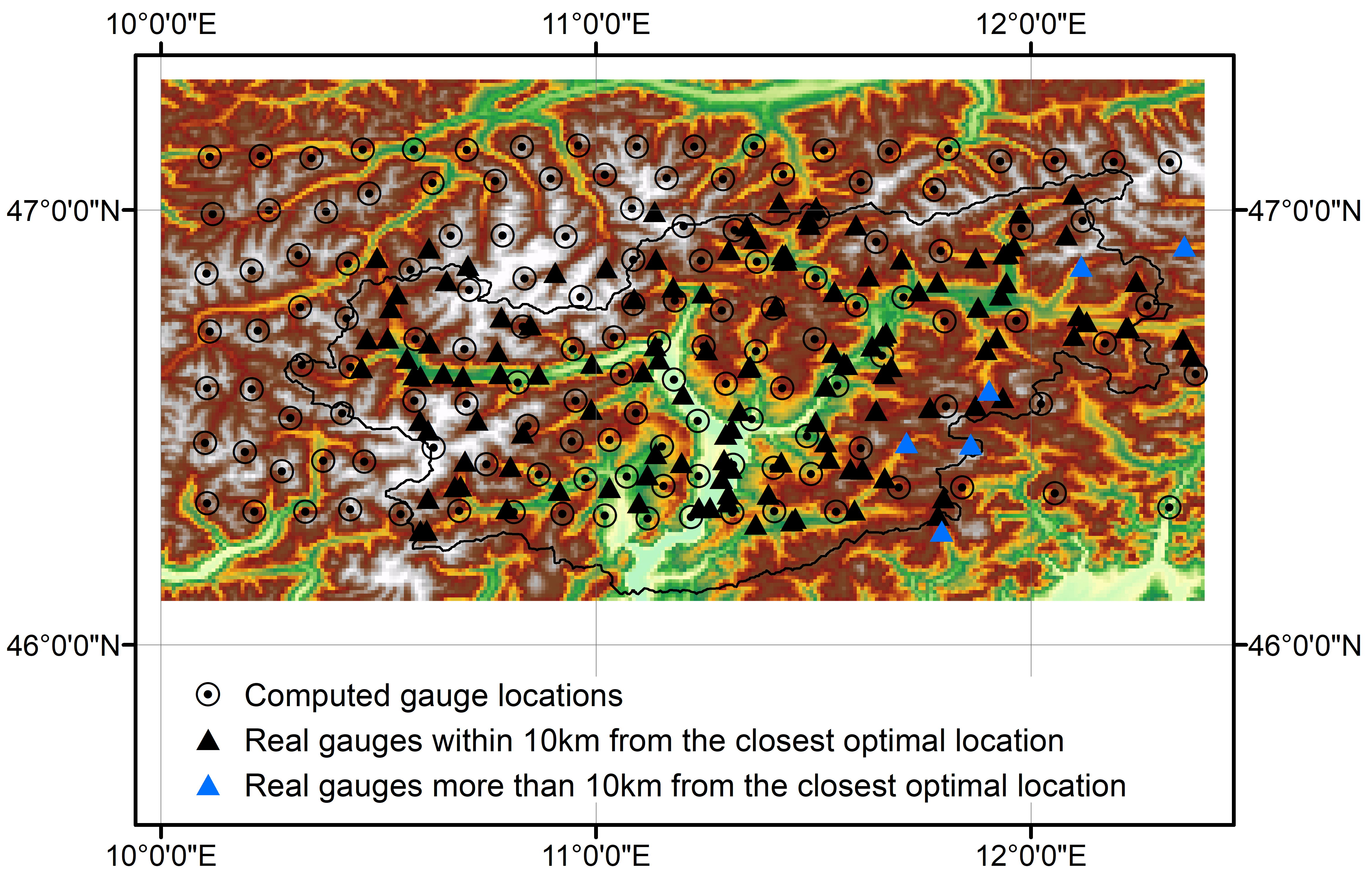}
\caption{Comparing existing locations with optimal placement for Adige region.}
\label{fig:compare_Adige}
\end{figure*}

In order to investigate how far the actual gauges are from their optimal locations, we look at Euclidean distance between a certain gauge and the closest optimal location. We pick a measure of closeness as a circle of a certain radius and count the number of gauges that fall within that distance from any optimal location. Tables \ref{tab1} and \ref{tab2} give direct counts of the gauges that are close and far from optimal locations, and Figures~\ref{fig:compare_OK} and \ref{fig:compare_Adige} provide the visual representation of the these results. 

Obviously, some locations may be inaccessible (especially in complex terrain) although optimal. The methodology introduced in this work can be modified to account for physical constraints and other resource allocation problems, but this is beyond the scope of this article and recommended for future work.

\begin{table}[h!]
\begin{center}
\begin{tabular}{ |c|c|c| } 
 \hline
 Radius (km)     &  gauges within radius $r$ &  gauges not within radius $d$ \\ 
 \hline
 $r=2$  & 8 & 130 \\
 $r=5$  & 60 & 78 \\ 
 $r=10$ & 132 & 6 \\ 
 \hline
\end{tabular}
\end{center}
\caption{Adige region: number of gauges close to optimal locations for different closeness measures.}
\label{tab2}
\end{table}

\begin{table}[h!]
\begin{center}
\begin{tabular}{ |c|c|c| } 
 \hline
 Radius (km)     & gauges within radius $r$ &  gauges not within radius $d$ \\ 
 \hline
 $r=5$  & 7 & 111 \\ 
 $r=10$ & 18 & 100 \\ 
  $r=15$ & 49 & 69 \\ 
  $r=15$ & 86 & 32 \\ 
 \hline
\end{tabular}
\end{center}
\caption{Oklahoma region: number of gauges close to optimal locations for different closeness measures}
\label{tab1}
\end{table}

\begin{figure}[ht]
\centering
\includegraphics[width=0.48\linewidth]{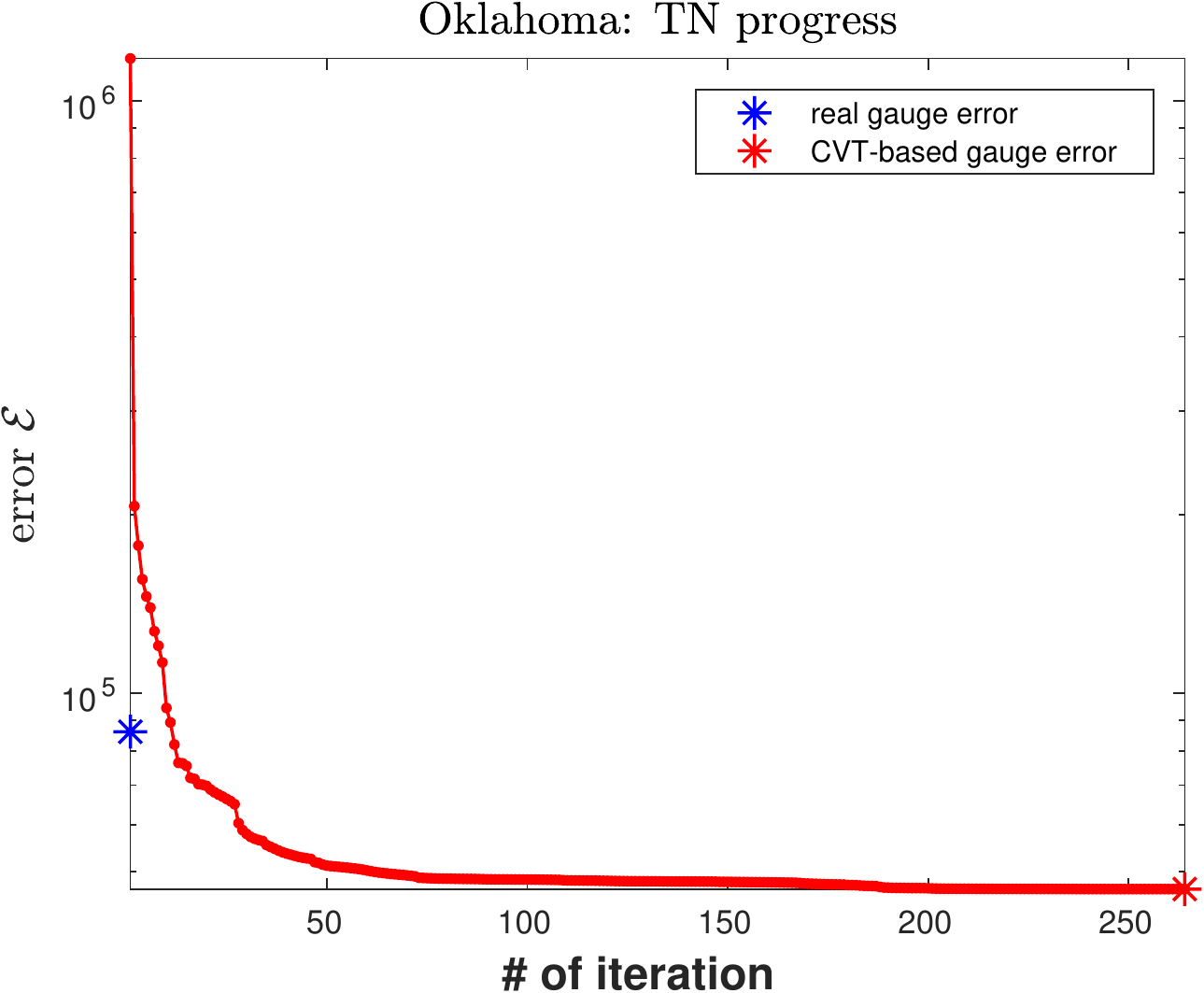}
\hfill
\includegraphics[width=0.48\linewidth]{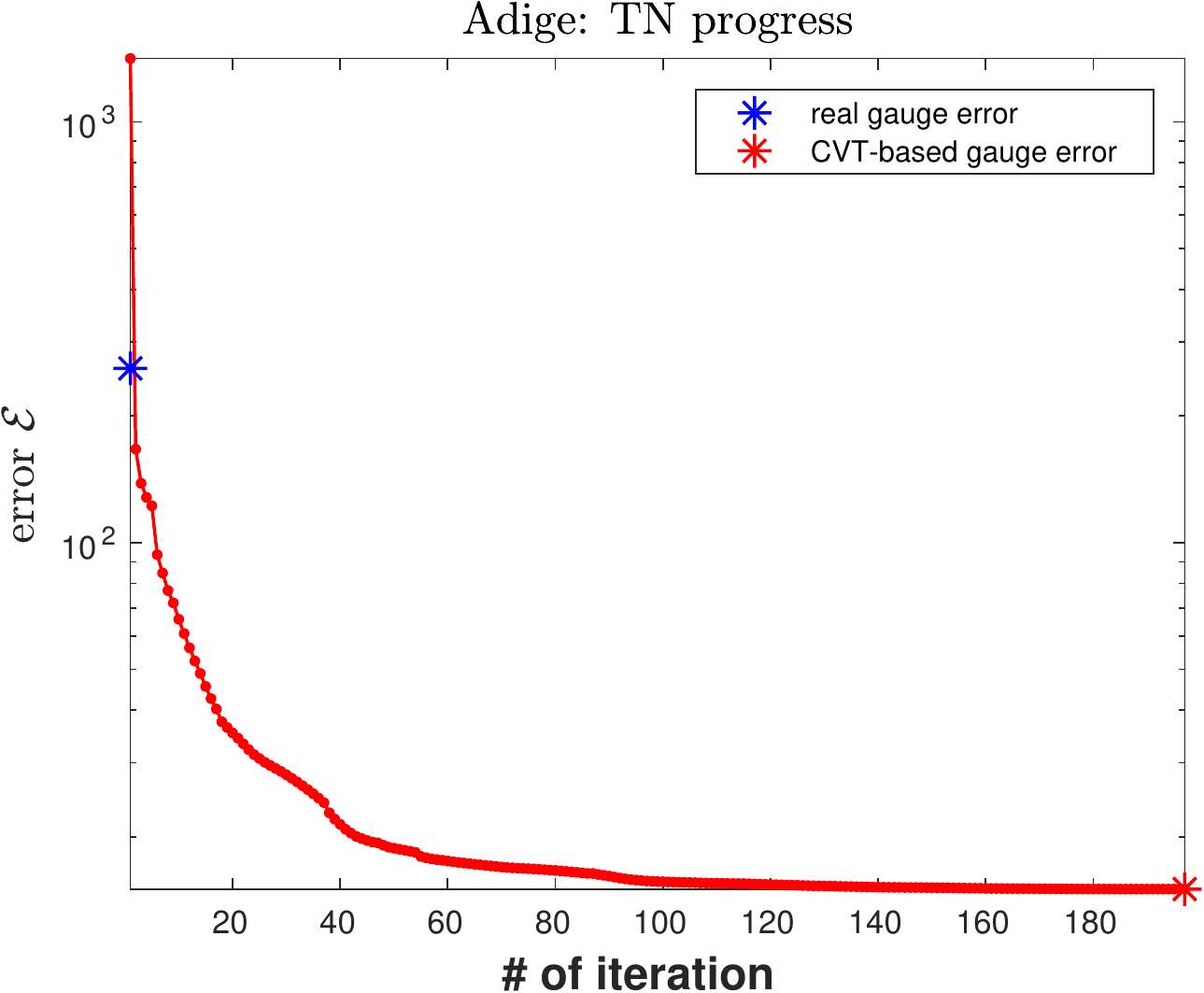}
\caption{Decay of the approximation error $\mathcal{E}$ \eqref{CVT_energy} during TN iteration for Oklahoma region (left) and Adige (right). CVT-based location denoted as a red asterisk gives a smaller error value comparing to the error computed for the real-gauge locations (shown as a blue asterisk).}
\label{fig: tn}
\end{figure}

As a measure of how well selected locations approximate the precipitation data in a given region, we can compute the error functional given by \eqref{CVT_energy}. As shown in Figure~\ref{fig: tn}, the approximation error is significantly decreased by running the optimization routine. More specifically, the error decreases from $259.56$ to $15.02$ for Adige data and from $8.61\cdot 10^4$ to $4.67\cdot 10^4$ for Oklahoma. Iteration was terminated once local minimum was achieved according to the chosen stopping criterion.

\section{Summary}
We developed an automated strategy that allows to find optimal precipitation gauge locations in any given region based on the variability pattern in the precipitation over the region. The proposed algorithm could potentially be applied to any precipitation dataset (including re-analysis products) that is long enough to capture the temporal variability of precipitation, regardless of any seasonality. We chose satellite-based observations because of their global (or quasi-global) coverage, making this method applicable anywhere else in the world. This is particularly useful when planning a field campaign to select sampling sites or when installing a new gauge network to pick the optimal number of gauges and their locations.

\section{Acknowledgements}
The authors are grateful to Ufficio Idrografico di Bolzano for providing Alto-Adige data. ME acknowledges support provided by the US National Science Foundation CAREER grant DMS-1056821 that funded Mason Modeling Days workshop where initial conceptualization of this work occurred. 



\bibliographystyle{elsarticle-num}
\bibliography{references}


\section{Appendix A: Optimization problem formulation.}
\label{app: opt}

Following \cite{Okabe}, let $z(\bx)$ be the random variable representing rainfall at location $\bx\in W$. We can think of it as 
$$
z(\bx)=m(\bx)+\epsilon(\bx)
$$
with $E(\epsilon(\bx))]=0$, so that $m(\bx)=E(z(\bx))$ represents the average trend of the process over the given region. 

Consider $k$ rain gauges (pluviometers) $\bx_1,\ldots, \bx_k$ placed within the region $W$ and let $V_i$ be the Voronoi (Thiessen) region associated with the $i$-th gauge.  As done by Thiessen back in 1911 \cite{thiessen}, we can estimate the total precipitation in the region $W$ as
$$
Z = \int_S{z(\bx)d\bx} = \sum^k_{i=1}\int_{V_i}{z(\bx)d\bx} \approx \sum^n_{i=1}|V_i| z(\bx_i)
$$

This leads to the natural way to formulate the optimization problem for optimal gauge placement to minimize the expected squared approximation error:
\begin{multline}
F(\bx_1,\ldots, \bx_k) = E[\sum^k_{i=1}\int_{V_i}(z(\bx)-z(\bx_i))^2 d\bx]  \\
=\sum^k_{i=1}\int_{V_i} [(m(\bx)-m(\bx_i))]^2 + [Var(\bx)-Var(\bx_i)]^2  \\
+ 2 Var(\bx) Var(\bx_i) [1-Corr(\bx,\bx_i)]\, d\bx \rightarrow \min
\label{full_opt}
\end{multline}
Here $Var(\bx)=E(\epsilon(\bx)^2)$ and $Var(\bx_i)=E(\epsilon(\bx_i)^2)$ denote variances at points $\bx$ and $\bx_i$ respectively and $Corr(\bx,\bx_i)=E[\epsilon(\bx)\epsilon(\bx_i)]/[Var(\bx)Var(\bx_i)]$ is the corresponding correlation of the time series at locations $\bx$ and $\bx_i$.

If the variance $Var(\bx)$ is small compared to the change in the average trend $m(\bx)$:
$$
F(\bx_1,\ldots, \bx_k) \approx  \sum^k_{i=1}\int_{V_i} [(m(\bx)-m(\bx_i))]^2\, d\bx \rightarrow \min
\label{opt1}
$$
One can solve this approximated problem to obtain optimal gauge placement in this case.

If on the other hand the change in the average trend $m(\bx)$ is small compared with the variance $Var(\bx)$, which is considered to be constant $Var(\bx)=\varepsilon$, the following approximation is valid:
$$
F(\bx_1,\ldots, \bx_k) \approx \sum^k_{i=1}\int_{V_i} {2 \varepsilon^2 [1-Corr(\bx,\bx_i)]\, d\bx} \rightarrow \min
\label{opt2}
$$

Hence, if there is a relatively small change in the means, we can solve this approximated problem instead.

Notice that in this case $F(\bx)$ is a variant of the CVT energy \eqref{CVT_energy1} 
with $f(||\bx-\bx_i||)=\varepsilon^2(1-Corr(\bx,\bx_i))$ and $\rho(\bx)=1$.

\section{Appendix B: Truncated-Newton Algorithm for CVT calculation}
\label{app:tn}
 
 Here, we give a brief review of  the truncated Newton algorithm. For more details we refer readers to \cite{nash2000survey}. To optimize a problem of the form
 $$
 \min_\bx f(\bx)
 $$
 at the $j$-th TN iteration a search direction $p$ is computed as an approximate solution to the Newton equations
 $$
 \nabla^2 f(\bx^{(j)}) p = -\nabla f(\bx^{(j)})
 $$
 where $\bx^{(j)}$ is the current approximation to the solution of the optimization problem.  The search direction $p$ is computed using the linear conjugate-gradient algorithm (CG).  The necessary Hessian-vector products are estimated using finite differencing.  The TN algorithm only requires that values of $f(\bx)$ and $\nabla f(\bx)$ are computed.  TN has low storage requirements, and has low computational costs per iteration, and hence is suitable for solving large optimization problems \cite{NaNo91a}.
 
The following are the steps necessary to compute discrete CVT using TN method using a pre-defined density function $\rho$:
 \begin{enumerate} 
 \renewcommand{\theenumi}{\roman{enumi}}
 \item Given the discrete set of points ${\bf Y}=\{y_i\}_{i=1}^{m}\in W$,
 \item Give the discrete energy function $$\cg ({\bx})=\sum_{ V_i} \sum_{ j}\rho (y_j)\|\bx_i -y_j\|^2$$ where $j$ is the index for those $y$ included in the voronoi set $V_i$,
 \item Compute its corresponding gradient value $$\nabla_i \cg ({\bx})=\sum_{ j}\rho (y_j)2(\bx_i -y_j),$$
 \item By Taylor series: $$\nabla \cg(\bx+\alpha v)=\nabla \cg(\bx)+\alpha \nabla^2 \cg(\bx) v$$ we can approximate the necessary component $\nabla^2 \cg(\bx) v$ of CG as $\ds\nabla^2 \cg(\bx) v=\frac{\nabla \cg(\bx+\alpha p)-\nabla \cg(\bx)}{\alpha}$,
 \item Substitue the above information to CG described as following to compute the search direction $p$:
 \begin{itemize}
 \item  $r_0:=-\nabla \cg(\bx)$
 \item $v_0:=r_0$
 \item $k:=0$
 \item repeat
 \begin{itemize}
 \item $\alpha_k := \frac{r_k^{T}r_k}{v_k^{T}\nabla^2 \cg(\bx) v_k}$
 \item $p_{k+1}:=p_k + \alpha_k v_k$
 \item $r_{k+1}:=r_k -\alpha_k \nabla^2 \cg(\bx) v_k$
 \item if $r_{k+1}$ is sufficiently small then exit loop, end if
 \item $\beta_k := \frac{r_{k+1}^{T}r_{k+1}}{r_k^{T}r_k}$
 \item $v_{k+1}:= r_{k+1}+\beta_k v_k$
 \item $k:=k+1$
 \end{itemize}
 \noindent end repeat
 \end{itemize}
 \item Test whether $\nabla \cg(\bx) p<0$. If so, accept $p$ as a descent direction, otherwise, take $p=-\nabla \cg(\bx)$,
 \item Use Armijo line search to determine the step size $\alpha$, then update $\bx$ by $\bx=\bx+\alpha p$,
 \item Go back to step iii) until a stopping criterion is reached.
 \end{enumerate}
 
\end{document}